\documentclass{article}

\usepackage{amsfonts}
\usepackage[all]{xy}
\usepackage{graphicx}
\usepackage{amsthm}
\usepackage{amsmath,amssymb,bbm}
\usepackage[colorlinks=true,citecolor=black,linkcolor=black,urlcolor=blue]{hyperref}
\usepackage{verbatim,fullpage,multicol}
\usepackage[table]{xcolor}
\usepackage[T1]{fontenc}
\usepackage{marvosym}
\usepackage{appendix}
\usepackage[utf8]{inputenc}
\usepackage{braket}
\usepackage{multirow}
\usepackage{rotating}
\usepackage{caption}
\usepackage{authblk}
\usepackage{subcaption}
\usepackage{ mathrsfs }
\usepackage{xcolor}
\usepackage{mathrsfs}
\usepackage{enumerate}
\usepackage[mathscr]{euscript}
\usepackage{float}
\usepackage{tikz}
\usepackage{tkz-graph}
\usepackage{bbold}

\newcommand{\Inc}{\mathrm{Inc}}

\definecolor{lightgray}{gray}{0.7}
\definecolor{midgray}{gray}{.9}

\tikzstyle{vertex}=[circle, draw, inner sep=0pt, minimum size=3pt]
\newcommand{\vertex}{\node[vertex]}
\tikzstyle{vtx}=[circle, draw, inner sep=0pt, minimum size=12pt]

\definecolor{darkgreen}{cmyk}{.9,0,.9,.2}

\newtheorem{theorem}{Theorem}[section]
\newtheorem{fact}[theorem]{Fact}
\newtheorem{proposition}[theorem]{Proposition}
\newtheorem{example}[theorem]{Example}
\newtheorem{notation}[theorem]{Notation}

\DeclareMathAlphabet{\mathbbold}{U}{bbold}{m}{n}

\def\P{{\mathbb P}}

\def\S{{{\mathbb{S}}}}
\def\Skn{{{\mathbb{S}}^{k}_{n}}}

\everymath{\displaystyle}

\newcommand{\ct}{\operatorname{Crit}}
\newcommand{\G}{\mathbf{G}^k_n}

\begin{document}

\title{Block Circulant Graphs and the \\Graphs of Critical Pairs of a Crown}

 \author[1]{Rebecca E. Garcia\thanks{{\href{mailto:mth\_reg@shsu.edu}{mth\_reg@shsu.edu}}}}
 \author[2]{Pamela E. Harris\thanks{{\href{mailto:peh2@williams.edu}{peh2@williams.edu}}}}
 \author[3]{Bethany Kubik\thanks{{\href{mailto:bakubik@d.umn.edu}{bakubik@d.umn.edu}}}}
 \author[4]{Joseph M. Pedersen\thanks{{\href{mailto:joseph.m.pedersen.mil@mail.mil}{joseph.m.pedersen.mil@mail.mil}}}} 
 \author[5]{Shannon Talbott\thanks{\textcolor{blue}{\href{mailto:talbotts@moravian.edu}{talbotts@moravian.edu}}}} 
 
\affil[1]{Department of Mathematics and Statistics, Sam Houston State University}
\affil[2]{Department of Mathematics and Statistics, Williams College}
\affil[3]{Mathematics and Statistics Department, University of Minnesota Duluth}
\affil[4]{Department of Mathematical Sciences, United States University}
\affil[5]{Mathematics and Computer Science Department, Moravian College}
 \renewcommand\Authands{ and }

\maketitle
\begin{abstract}
 In this paper, we provide a natural bijection between a special family of block circulant graphs and the graphs of critical pairs of the posets known as generalized crowns. In particular, every graph in this family of block circulant graphs we investigate has a generating block row that follows a symmetric growth pattern of the all ones matrix. The natural bijection provides an upper bound on the chromatic number for this infinite family of graphs. 
\end{abstract}



\section{Introduction}\label{intro}
Circulant and block circulant matrices are well-studied families of matrices with applications to a variety of mathematical areas. In numerical analysis, circulant matrices play an important role as they can be diagonalized via a discrete Fourier transform \cite{Bini,Gray}. In cryptography, circulant matrices are used in the MixColumns step of the Advanced Encryption Standard and in error correcting codes, \cite{Daemen,Milenkovic}, while block circulant matrices have applications in Quasi-Cyclic (QC) codes \cite{Tanner}. In terms of graph theory,
Codenotti, Gerace, and Vigna showed that computing the chromatic number of circulant graphs, graphs whose adjacency matrices are circulant, is an NP-hard problem \cite{Codenotti}. This led to the development of efficient algorithms to compute the chromatic number of circulant graphs which improve current graph coloring algorithms \cite{Discepoli}. However, even with this improvement, the chromatic number problem for circulant and block circulant matrices remains an active area of research. 

In this paper, we introduce a special family of block circulant graphs, whose elements are denoted by $\mathcal{BC}_s^t$ and are parameterized by two positive integer values $s$ and $t$. These graphs are generated by a block row with a symmetric growth pattern of the all ones matrix. We then provide a connection between this family of block circulant graphs and a particular family of height 2 posets, the elements of which are called generalized crowns and are denoted by $\Skn$. 
Our main result is as follows:

\begin{theorem}\label{main}[Graph Isomorphism] Let $n \geq 3$ and $k \geq 0$. If $\G$ is the graph of critical pairs of $\Skn$, then $\G$ is graph isomorphic to
$\mathcal{BC}_{n+k}^{k+1}$.
\end{theorem}

The significance of this work comes from the relation between the chromatic number of the graph $\G$ and the order dimension of the crown $\Skn$. Felsner and Trotter showed that for every finite poset $\P$, 
$\dim(\P)\geq \chi(\mathbf{G}_{\P}^c)$ where $\dim(\P)$ denotes the order dimension of the poset $\P$, $\mathbf{G}_{\P}^c$ denotes the graph of critical pairs of $\P$, and $\chi$ refers to the chromatic number of $\mathbf{G}_{\P}^c$ \cite[Lemma~3.3]{FeTr00}.
This inequality produces an upper bound for the chromatic number for all block circulant graphs $\mathcal{BC}_s^t$.

This paper is organized as follows: Section \ref{sec1} provides the necessary background material on block circulant matrices and poset theory to make our approach precise. In Section \ref{sec2}, we prove the graph isomorphism between $\G$ and $\mathcal{BC}_{n+k}^{k+1}$ (see Theorem \ref{main}). Section \ref{sec3} concludes with a few open questions and directions for future work.  

\section{Background}\label{sec1}
An \emph{$m$-block circulant matrix} $\mathcal{C}$ is a matrix of dimension $nm\times nm$ that is generated by the matrices $C_1,C_2,\ldots,C_n$ of dimension $m\times m$, where the block rows of $\mathcal{C}$ are obtained by cyclically shifting the $C_i$'s as follows:
\[\mathcal{C}=\text{circ}(C_1,C_2,\dots,C_n)=\begin{bmatrix}
C_{1} & C_{2}&C_{3}&\cdots&C_{n-1}&C_{n}\\
C_{n} & C_{1}&C_{2}&\cdots&C_{n-2}&C_{n-1}\\
\vdots&\vdots&\vdots&\ddots&\vdots&\vdots\\
C_{3}&C_{4}&C_{5}&\cdots&C_{1}&C_{2}\\
C_{2}&C_{3}&C_{4}&\cdots&C_n&C_{1}
\end{bmatrix}.\]
The matrix $\mathcal{C}$ has \emph{block size} $m$ and a \emph{generating block row} consisting of the matrices $C_1,C_2,\ldots,C_n$. 
If the matrices $C_1,C_2,\ldots,C_n$ are circulant, then $C$ is said to be an \emph{$m$-block circulant matrix with circulant blocks}. 
We note that a circulant matrix is a $1$-block circulant matrix. Whenever $m$ is understood, we refer to an $m$-block circulant matrix as a block circulant matrix. 

Let $\mathcal{BC}$ denote the family of block circulant matrices with non-negative integer entries. Abusing notation, let $\mathcal{BC}$ also denote the family of graphs whose adjacency matrices are block circulant. 
Our object of study is an infinite subfamily of graphs in $\mathcal{BC}$, whose elements are the graphs denoted by $\mathcal{BC}_s^t$ having a $t$-block circulant adjacency matrix and generating block row $B_1,B_2,\ldots, B_s$, where $t\geq 1$ and $s\geq t+2$. To define $B_1,B_2,\ldots, B_s$, set the following notation:

\begin{notation} Let $1\leq i\leq t$. 
The $t\times t$ 
block  $\mathbbold{1}^{i}$ has an $i \times i$ block 
of ones in the upper right corner, with the remaining entries of the $t \times t$ block being zero.  
Similarly, the $t \times t$ block $_{i}\mathbbold{1}$ has an $i \times i$ block of ones in the lower left corner, with the remaining entries of the $t \times t$ block being zero. For example, a $6 \times 6$ block $\mathbbold{1}^{4}$ and a $6 \times 6$ block $_{2}\mathbbold{1}$ are shown below:\\
\begin{minipage}[h]{\textwidth}
\begin{minipage}{0.5\textwidth}
\[ \mathbbold{1}^{4} = \begin{bmatrix} 0 & 0 & 1 & 1 & 1 & 1 \\  0 & 0 & 1 & 1 & 1 & 1 \\  0 & 0 & 1 & 1 & 1 & 1 \\  0 & 0 & 1 & 1 & 1 & 1 \\ 0 & 0 & 0 & 0 & 0 & 0 \\ 0 & 0 & 0 & 0 & 0 & 0   \end{bmatrix} \]
\end{minipage}
\hspace{-1in}
\begin{minipage}{0.5\textwidth}
\[ _{2}\mathbbold{1} = \begin{bmatrix}  0 & 0 & 0 & 0 & 0 & 0 \\ 0 & 0 & 0 & 0 & 0 & 0  \\ 0 & 0 & 0 & 0 & 0 & 0 \\ 0 & 0 & 0 & 0 & 0 & 0 \\ 1 & 1 & 0 & 0 & 0 & 0  \\  1 & 1 & 0 & 0 & 0 & 0 \end{bmatrix}. \]
\end{minipage}
\end{minipage}
\noindent
The $t \times t$ block $\mathbbold{1}_{t\times t}$ is the all ones $t \times t$ matrix. 
Similarly, the $t \times t$ block $\mathbbold{0}_{t\times t}$ is the zero $t \times t$ matrix. 
\end{notation}

Define the generating block row of $\mathcal{BC}_s^t$ as follows:

\noindent
{\bf Case 1:} If $s\geq 2t$, then set
\begin{align*}
B_i&=\begin{cases}\mathbbold{0}_{t\times t} & \mbox{if  $i=1$}\\
\mathbbold{1}^{i-1}&\mbox{if $2\leq i\leq t$}\\
\mathbbold{1}_{t\times t}&\mbox{if $t+1\leq i\leq s-t+1$}\\
_{s+1-i}\mathbbold{1}&\mbox{if $s-t+2\leq i\leq s$}.
\end{cases}
\end{align*}
As an example, the generating block row of $\mathcal{BC}_7^3$ is 
\[
\left[\arraycolsep=2pt\def\arraystretch{1}
\begin{array}{c|c|c|c|c|c|c}
B_1&B_2&B_3&B_4&B_5&B_6&B_7
\end{array}\right]=
\left[\arraycolsep=1.4pt\def\arraystretch{1}
\begin{array}{llp{3mm}|@{\hskip 2mm}ccp{3mm}|@{\hskip 2mm}ccp{3mm}|@{\hskip 2mm}ccp{3mm}|@{\hskip 2mm}ccp{3mm}|@{\hskip 2mm}ccp{3mm}|@{\hskip 2mm}ccc}
0&0&0&0&0&1&0&1&1&1&1&1&1&1&1&0&0&0&0&0&0\\
0&0&0&0&0&0&0&1&1&1&1&1&1&1&1&1&1&0&0&0&0\\
0&0&0&0&0&0&0&0&0&1&1&1&1&1&1&1&1&0&1&0&0
\end{array}
\right].\]
\noindent
{\bf Case 2:} If $s\leq 2t-1$, then set 
\begin{align*}
B_i&=\begin{cases}\mathbbold{0}_{t\times t} & \mbox{if  $i=1$}\\
\mathbbold{1}^{i-1}&\mbox{if $2\leq i\leq s-t+1$}\\
\begin{bmatrix}0&0&\mathbbold{0}_{(t-s-1+i)\times (t-s-1+i)}\\0&\mathbbold{1}_{(s-t)\times (s-t)}&0\\
\mathbbold{0}_{(t-i+1)\times (t-i+1)}&0&0\end{bmatrix}
&\mbox{if $s-t+2\leq i\leq t$}\\
_{s-i+1}\mathbbold{1}&\mbox{if $t+1\leq i\leq s$}.
\end{cases}
\end{align*}
\noindent
As an example, the generating block row of $\mathcal{BC}_6^4$ is 
\[
\left[\arraycolsep=2.0pt\def\arraystretch{1.2}
\begin{array}{c|c|c|c|c|c}
B_1&B_2&B_3&B_4&B_5&B_6
\end{array}\right]=
\left[\arraycolsep=1.4pt\def\arraystretch{1.2}
\begin{array}{lllp{3mm}|@{\hskip 2mm}cccp{3mm}|@{\hskip 2mm}cccp{3mm}|@{\hskip 2mm}cccp{3mm}|@{\hskip 2mm}cccp{3mm}|@{\hskip 2mm}cccc}
0&0&0&0&
0&0&0&1&
0&0&1&1&
0&0&0&0&
0&0&0&0&
0&0&0&0\\
0&0&0&0&
0&0&0&0&
0&0&1&1&
0&1&1&0&
0&0&0&0&
0&0&0&0\\
0&0&0&0&
0&0&0&0&
0&0&0&0&
0&1&1&0&
1&1&0&0&
0&0&0&0\\
0&0&0&0&
0&0&0&0&
0&0&0&0&
0&0&0&0&
1&1&0&0&
1&0&0&0
\end{array}
\right].\]

To provide the connection to the graph of critical pairs of generalized crowns, we now focus our attention on the necessary background in poset theory.
We assume some familiarity with posets, order dimension, and chromatic number and refer the interested reader to \cite{Diestel,Tr92,Tr95} for further background. 

Throughout this paper, let $n,k\in \mathbb{N}$ with $n\geq 3$ and $k\geq 0$.  
The \emph{generalized crown}, denoted $\Skn$, is a height 2 poset with $\min(\Skn)=\{a_1,\ldots,a_{n+k}\}$ and $\max(\Skn)=\{b_1,\ldots,b_{n+k}\}$, where 
\begin{enumerate}[\rm(1)]
\item $b_i||a_i,a_{i+1},\ldots,a_{i+k}$, and
\item $b_i>a_{i+k+1},a_{i+k+2}, \ldots,a_{i-1}$.
\end{enumerate}

\begin{figure}[H]
\centering
 \begin{tikzpicture}
 \vertex[fill,label=below:$a_1$](a1) at (1,0) {};
 \vertex[fill,label=below:$a_2$](a2) at (2,0) {};
 \vertex[fill,label=below:$a_3$](a3) at (3,0) {};
 \vertex[fill,label=below:$a_4$](a4) at (4,0) {};
 \vertex[fill,label=below:$a_5$](a5) at (5,0) {};
  \vertex[fill,label=below:$a_6$](a6) at (6,0) {};
   \vertex[fill,label=below:$a_7$](a7) at (7,0) {};
  \vertex[fill,label=above:$b_1$](b1) at (1,1) {};
 \vertex[fill,label=above:$b_2$](b2) at (2,1) {};
 \vertex[fill,label=above:$b_3$](b3) at (3,1) {};
 \vertex[fill,label=above:$b_4$](b4) at (4,1) {};
 \vertex[fill,label=above:$b_5$](b5) at (5,1) {};
 \vertex[fill,label=above:$b_6$](b6) at (6,1) {};
 \vertex[fill,label=above:$b_7$](b7) at (7,1) {};
	\draw (1,1)--(5,0) (1,1)--(6,0) (1,1)--(7,0);
	\draw (2,1)--(6,0) (2,1)--(7,0) (2,1)--(1,0);
	\draw (3,1)--(7,0) (3,1)--(1,0) (3,1)--(2,0);
  	\draw (4,1)--(1,0) (4,1)--(2,0) (4,1)--(3,0);
	\draw (5,1)--(2,0) (5,1)--(3,0) (5,1)--(4,0);
	\draw (6,1)--(3,0) (6,1)--(4,0) (6,1)--(5,0);
	\draw (7,1)--(4,0) (7,1)--(5,0) (7,1)--(6,0);
	 \end{tikzpicture}
\caption{The crown $\S_4^3$}\label{prettycrown}
\end{figure}
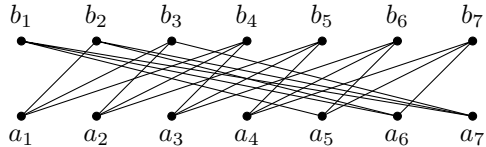

Figure \ref{prettycrown} provides the diagram for the generalized crown $\mathbb{S}^3_4$. Identifying $a_i$ with $a_{i-(n+k)}$ and $b_i$ with $b_{i-(n+k)}$ for $i>n+k$ is called \emph{cyclic indexing}. 
The set of all incomparable pairs of $\Skn$ is denoted by
$\Inc(\Skn)=\{(x,y)\in\Skn\times \Skn\colon x||y\}.$
The pair $(x,y)\in \Skn\times\Skn$ is \emph{critical} if the following conditions hold:(i) $x||y$;
(ii) $D(x)\subset D(y)$; and
(iii) $U(y)\subset U(x)$,
where $D(u)=\{z\in\Skn:z<u\}$ and $U(w)=\{z\in\Skn:w<z\}$ for any $u,w\in\Skn$. Let $\ct(\Skn)$ denote the set of all critical pairs of $\Skn$. 
An \emph{alternating cycle} is a sequence $\{(x_i,y_i) : 1\leq i\leq k\}$
of ordered pairs from $\Inc(\Skn)$, where $y_i\leq x_{i+1}$ in $\Skn$ (cyclically) for
$i=1,2,\dots,k$. An alternating cycle is said to be \emph{strict} if $y_i\leq x_j$ in $\Skn$ if and only if $j=i+1$ (cyclically) for
$i,j=1,2,\dots, k.$

The \emph{strict hypergraph of critical pairs of $\Skn$},
denoted $\mathbf{H}_{n}^k$, is the hypergraph with vertices $\ct(\Skn)$ and edges consisting of subsets of $\ct(\Skn)$
whose duals form strict alternating cycles. If $(x,y)$ is a critical pair, then $(y,x)$ is its dual. Let $\mathbf{G}_{n}^k$ denote the graph of $\mathbf{H}_{n}^k$. That is, $\mathbf{G}_n^k$ is a graph with vertices $\ct(\Skn)$ and edges consisting of size 2 subsets of $\ct(\Skn)$
whose duals form strict alternating cycles. With these definitions at hand, we now formulate the set bijection between the families of graphs $\{\G\}_{n\geq 3, k\geq 0}$ and $\{\mathcal{BC}_{s}^t\}_{t\geq 1,s\geq t+2}$.

\section{The graph isomorphism of $\G$ and $\mathcal{BC}_s^t$}
\label{sec2}

In this section, we prove: Propositions \ref{blockcirculant1} and \ref{blockcirculant2}, which make clear that all $\G$ belongs to the family $\{\mathcal{BC}_{s}^t\}$, and
Theorem \ref{blockcirc}, which demonstrates that every graph $\mathcal{BC}_s^t$ arises as the graph $\mathbf{G}_{t-1}^{s-t+1}$, where $t\geq1$ and $s\geq t+2$. Together these results establish Theorem~\ref{main}.

Let $\mathcal{A}_n^k$ denote the adjacency matrix of $\G$. To give the entries of the matrix $\mathcal{A}_n^k$, first note that $\Skn$ has  $(n+k)(k+1)$ critical pairs, 
which we list in lexicographical order on their dual and use this labeling on the rows (and by symmetry columns) of the matrix $\mathcal{A}_n^k$. Our notation is as follows:

\begin{notation}\label{matrix-index}
Fix $1\leq i,j\leq n+k$ and let $A_{i,j}$ denote the $(k+1)\times(k+1)$ submatrix whose rows are labeled by the $k+1$ critical pairs: $(a_{i},b_i),(a_{i+1},b_i),\ldots,(a_{i+k-1},b_i),(a_{i+k},b_i),$ and  whose columns are labeled by the $k+1$ critical pairs: $(a_j,b_j),(a_{j+1},b_j),\ldots,(a_{j+k-1},b_j),(a_{j+k},b_j),$ where all the subscripts of the first component are taken cyclically modulo $n+k$.
Then $\mathcal{A}_n^k=[A_{i,j}]_{1\leq i,j\leq n+k}.$
For fixed $1\leq i,j\leq n+k$, let $m_{u,v}$ denote the $(u,v)$-entry of the submatrix $A_{i,j}$. Notice that $u$ ranges from $i$ to $k+i$, where the order is fixed and all terms are taken modulo $n+k$.
Similarly, $v$ ranges from $j$ to $j+k$, where the order is fixed and the terms are taken modulo $n+k$. Denote these ranges by writing 
$u\in [i, i+1,\ldots, i+k]_{\!\!\mod{(n+k)}}$ and $v\in[j, j+1,\ldots, j+k]_{\!\!\mod{(n+k)}}.$
\end{notation}

\begin{example}
The matrix $\mathcal{A}_4^3$ is determined by which duals of critical pairs of $\S^3_4$ form strict alternating cycles of size 2; see Table \ref{prettytable}. 
A computation shows that $\mathbf{G}^3_4$ is 3-colorable; see Figure \ref{prettygraph}.

\noindent
\begin{minipage}[b]{\textwidth} 
\noindent
\begin{minipage}[b]{0.57\textwidth}
\begin{table}[H]
\centering
 \resizebox{\textwidth}{!}{
\begin{tabular}{c||p{1mm}p{1mm}p{1mm}p{1mm}@{\hskip 2.5mm}|p{1mm}p{1mm}p{1mm}p{1mm}@{\hskip 2.5mm}|p{1mm}p{1mm}p{1mm}p{1mm}@{\hskip 2.5mm}|p{1mm}p{1mm}p{1mm}p{1mm}@{\hskip 2.5mm}|p{1mm}p{1mm}p{1mm}p{1mm}@{\hskip 2.5mm}|p{1mm}p{1mm}p{1mm}p{1mm}@{\hskip 2.5mm}|p{1mm}p{1mm}p{1mm}p{1mm}@{\hskip 2.5mm}|}
\
&	\begin{turn}{-90}$\hspace{-2mm}(a_1,b_1)\hspace{.75mm} $\end{turn}	&	\begin{turn}{-90}$\hspace{-2mm}(a_2,b_1)\hspace{.75mm} $\end{turn}	&	\begin{turn}{-90}$\hspace{-2mm}(a_3,b_1)\hspace{.75mm} $\end{turn}	&	\begin{turn}{-90}$\hspace{-2mm}(a_4,b_1)\hspace{.75mm} $\end{turn}	&	\begin{turn}{-90}$\hspace{-2mm}(a_2,b_2)\hspace{.75mm} $\end{turn}	&	\begin{turn}{-90}$\hspace{-2mm}(a_3,b_2)\hspace{.75mm} $\end{turn}	&	\begin{turn}{-90}$\hspace{-2mm}(a_4,b_2)\hspace{.75mm} $\end{turn}	&	\begin{turn}{-90}$\hspace{-2mm}(a_5,b_2)\hspace{.75mm} $\end{turn}	&	\begin{turn}{-90}$\hspace{-2mm}(a_3,b_3)\hspace{.75mm} $\end{turn}	&	\begin{turn}{-90}$\hspace{-2mm}(a_4,b_3)\hspace{.75mm} $\end{turn}	&	\begin{turn}{-90}$\hspace{-2mm}(a_5,b_3)\hspace{.75mm} $\end{turn}	&	\begin{turn}{-90}$\hspace{-2mm}(a_6,b_3)\hspace{.75mm} $\end{turn}	&	\begin{turn}{-90}$\hspace{-2mm}(a_4,b_4)\hspace{.75mm} $\end{turn}	&	\begin{turn}{-90}$\hspace{-2mm}(a_5,b_4)\hspace{.75mm} $\end{turn}	&	\begin{turn}{-90}$\hspace{-2mm}(a_6,b_4)\hspace{.75mm} $\end{turn}	&	\begin{turn}{-90}$\hspace{-2mm}(a_7,b_4)\hspace{.75mm} $\end{turn}	&	\begin{turn}{-90}$\hspace{-2mm}(a_5,b_5)\hspace{.75mm} $\end{turn}	&	\begin{turn}{-90}$\hspace{-2mm}(a_6,b_5)\hspace{.75mm} $\end{turn}	&	\begin{turn}{-90}$\hspace{-2mm}(a_7,b_5)\hspace{.75mm} $\end{turn}	&	\begin{turn}{-90}$\hspace{-2mm}(a_1,b_5)\hspace{.75mm} $\end{turn}	&	\begin{turn}{-90}$\hspace{-2mm}(a_6,b_6)\hspace{.75mm} $\end{turn}	&	\begin{turn}{-90}$\hspace{-2mm}(a_7,b_6)\hspace{.75mm} $\end{turn}	&	\begin{turn}{-90}$\hspace{-2mm}(a_1,b_6)\hspace{.75mm} $\end{turn}	&	\begin{turn}{-90}$\hspace{-2mm}(a_2,b_6)\hspace{.75mm} $\end{turn}	&	\begin{turn}{-90}$\hspace{-2mm}(a_7,b_7)\hspace{.75mm} $\end{turn}	&	\begin{turn}{-90}$\hspace{-2mm}(a_1,b_7)\hspace{.75mm} $\end{turn}	&	\begin{turn}{-90}$\hspace{-2mm}(a_2,b_7)\hspace{.75mm} $\end{turn}	&	\begin{turn}{-90}$\hspace{-2mm}(a_3,b_7)\hspace{.75mm} $\end{turn}	\\\hline\hline
$(a_1,b_1)$	&	0	&	0	&	0	&	0	&	0	&	0	&	0	&	\cellcolor{lightgray}1	&	0	&	0	&	\cellcolor{lightgray}1	&	\cellcolor{lightgray}1	&	0	&	\cellcolor{lightgray}1	&	\cellcolor{lightgray}1	&	\cellcolor{lightgray}1	&	0	&	0	&	0	&	0	&	0	&	0	&	0	&	0	&	0	&	0	&	0	&	0	\\	
$(a_2,b_1)$	&	0	&	0	&	0	&	0	&	0	&	0	&	0	&	0	&	0	&	0	&	\cellcolor{lightgray}1	&	\cellcolor{lightgray}1	&	0	&	\cellcolor{lightgray}1	&	\cellcolor{lightgray}1	&	\cellcolor{lightgray}1	&	\cellcolor{lightgray}1	&	\cellcolor{lightgray}1	&	\cellcolor{lightgray}1	&	0	&	0	&	0	&	0	&	0	&	0	&	0	&	0	&	0	\\	
$(a_3,b_1)$	&	0	&	0	&	0	&	0	&	0	&	0	&	0	&	0	&	0	&	0	&	0	&	0	&	0	&	\cellcolor{lightgray}1	&	\cellcolor{lightgray}1	&	\cellcolor{lightgray}1	&	\cellcolor{lightgray}1	&	\cellcolor{lightgray}1	&	\cellcolor{lightgray}1	&	0	&	\cellcolor{lightgray}1	&	\cellcolor{lightgray}1	&	0	&	0	&	0	&	0	&	0	&	0	\\	
$(a_4,b_1)$	&	0	&	0	&	0	&	0	&	0	&	0	&	0	&	0	&	0	&	0	&	0	&	0	&	0	&	0	&	0	&	0	&	\cellcolor{lightgray}1	&	\cellcolor{lightgray}1	&	\cellcolor{lightgray}1	&	0	&	\cellcolor{lightgray}1	&	\cellcolor{lightgray}1	&	0	&	0	&	\cellcolor{lightgray}1	&	0	&	0	&	0	\\	\hline	
$(a_2,b_2)$	&	0	&	0	&	0	&	0	&	0	&	0	&	0	&	0	&	0	&	0	&	0	&	\cellcolor{lightgray}1	&	0	&	0	&	\cellcolor{lightgray}1	&	\cellcolor{lightgray}1	&	0	&	\cellcolor{lightgray}1	&	\cellcolor{lightgray}1	&	\cellcolor{lightgray}1	&	0	&	0	&	0	&	0	&	0	&	0	&	0	&	0	\\	
$(a_3,b_2)$	&	0	&	0	&	0	&	0	&	0	&	0	&	0	&	0	&	0	&	0	&	0	&	0	&	0	&	0	&	\cellcolor{lightgray}1	&	\cellcolor{lightgray}1	&	0	&	\cellcolor{lightgray}1	&	\cellcolor{lightgray}1	&	\cellcolor{lightgray}1	&	\cellcolor{lightgray}1	&	\cellcolor{lightgray}1	&	\cellcolor{lightgray}1	&	0	&	0	&	0	&	0	&	0	\\	
$(a_4,b_2)$	&	0	&	0	&	0	&	0	&	0	&	0	&	0	&	0	&	0	&	0	&	0	&	0	&	0	&	0	&	0	&	0	&	0	&	\cellcolor{lightgray}1	&	\cellcolor{lightgray}1	&	\cellcolor{lightgray}1	&	\cellcolor{lightgray}1	&	\cellcolor{lightgray}1	&	\cellcolor{lightgray}1	&	0	&	\cellcolor{lightgray}1	&	\cellcolor{lightgray}1	&	0	&	0	\\	
$(a_5,b_2)$	&	\cellcolor{lightgray}1	&	0	&	0	&	0	&	0	&	0	&	0	&	0	&	0	&	0	&	0	&	0	&	0	&	0	&	0	&	0	&	0	&	0	&	0	&	0	&	\cellcolor{lightgray}1	&	\cellcolor{lightgray}1	&	\cellcolor{lightgray}1	&	0	&	\cellcolor{lightgray}1	&	\cellcolor{lightgray}1	&	0	&	0	\\	\hline	
$(a_3,b_3)$	&	0	&	0	&	0	&	0	&	0	&	0	&	0	&	0	&	0	&	0	&	0	&	0	&	0	&	0	&	0	&	\cellcolor{lightgray}1	&	0	&	0	&	\cellcolor{lightgray}1	&	\cellcolor{lightgray}1	&	0	&	\cellcolor{lightgray}1	&	\cellcolor{lightgray}1	&	\cellcolor{lightgray}1	&	0	&	0	&	0	&	0	\\	
$(a_4,b_3)$	&	0	&	0	&	0	&	0	&	0	&	0	&	0	&	0	&	0	&	0	&	0	&	0	&	0	&	0	&	0	&	0	&	0	&	0	&	\cellcolor{lightgray}1	&	\cellcolor{lightgray}1	&	0	&	\cellcolor{lightgray}1	&	\cellcolor{lightgray}1	&	\cellcolor{lightgray}1	&	\cellcolor{lightgray}1	&	\cellcolor{lightgray}1	&	\cellcolor{lightgray}1	&	0	\\	
$(a_5,b_3)$	&	\cellcolor{lightgray}1	&	\cellcolor{lightgray}1	&	0	&	0	&	0	&	0	&	0	&	0	&	0	&	0	&	0	&	0	&	0	&	0	&	0	&	0	&	0	&	0	&	0	&	0	&	0	&	\cellcolor{lightgray}1	&	\cellcolor{lightgray}1	&	\cellcolor{lightgray}1	&	\cellcolor{lightgray}1	&	\cellcolor{lightgray}1	&	\cellcolor{lightgray}1	&	0	\\	
$(a_6,b_3)$	&	\cellcolor{lightgray}1	&	\cellcolor{lightgray}1	&	0	&	0	&	\cellcolor{lightgray}1	&	0	&	0	&	0	&	0	&	0	&	0	&	0	&	0	&	0	&	0	&	0	&	0	&	0	&	0	&	0	&	0	&	0	&	0	&	0	&	\cellcolor{lightgray}1	&	\cellcolor{lightgray}1	&	\cellcolor{lightgray}1	&	0	\\	\hline	
$(a_4,b_4)$	&	0	&	0	&	0	&	0	&	0	&	0	&	0	&	0	&	0	&	0	&	0	&	0	&	0	&	0	&	0	&	0	&	0	&	0	&	0	&	\cellcolor{lightgray}1	&	0	&	0	&	\cellcolor{lightgray}1	&	\cellcolor{lightgray}1	&	0	&	\cellcolor{lightgray}1	&	\cellcolor{lightgray}1	&	\cellcolor{lightgray}1	\\	
$(a_5,b_4)$	&	\cellcolor{lightgray}1	&	\cellcolor{lightgray}1	&	\cellcolor{lightgray}1	&	0	&	0	&	0	&	0	&	0	&	0	&	0	&	0	&	0	&	0	&	0	&	0	&	0	&	0	&	0	&	0	&	0	&	0	&	0	&	\cellcolor{lightgray}1	&	\cellcolor{lightgray}1	&	0	&	\cellcolor{lightgray}1	&	\cellcolor{lightgray}1	&	\cellcolor{lightgray}1	\\	
$(a_6,b_4)$	&	\cellcolor{lightgray}1	&	\cellcolor{lightgray}1	&	\cellcolor{lightgray}1	&	0	&	\cellcolor{lightgray}1	&	\cellcolor{lightgray}1	&	0	&	0	&	0	&	0	&	0	&	0	&	0	&	0	&	0	&	0	&	0	&	0	&	0	&	0	&	0	&	0	&	0	&	0	&	0	&	\cellcolor{lightgray}1	&	\cellcolor{lightgray}1	&	\cellcolor{lightgray}1	\\	
$(a_7,b_4)$	&	\cellcolor{lightgray}1	&	\cellcolor{lightgray}1	&	\cellcolor{lightgray}1	&	0	&	\cellcolor{lightgray}1	&	\cellcolor{lightgray}1	&	0	&	0	&	\cellcolor{lightgray}1	&	0	&	0	&	0	&	0	&	0	&	0	&	0	&	0	&	0	&	0	&	0	&	0	&	0	&	0	&	0	&	0	&	0	&	0	&	0	\\	\hline	
$(a_5,b_5)$	&	0	&	\cellcolor{lightgray}1	&	\cellcolor{lightgray}1	&	\cellcolor{lightgray}1	&	0	&	0	&	0	&	0	&	0	&	0	&	0	&	0	&	0	&	0	&	0	&	0	&	0	&	0	&	0	&	0	&	0	&	0	&	0	&	\cellcolor{lightgray}1	&	0	&	0	&	\cellcolor{lightgray}1	&	\cellcolor{lightgray}1	\\	
$(a_6,b_5)$	&	0	&	\cellcolor{lightgray}1	&	\cellcolor{lightgray}1	&	\cellcolor{lightgray}1	&	\cellcolor{lightgray}1	&	\cellcolor{lightgray}1	&	\cellcolor{lightgray}1	&	0	&	0	&	0	&	0	&	0	&	0	&	0	&	0	&	0	&	0	&	0	&	0	&	0	&	0	&	0	&	0	&	0	&	0	&	0	&	\cellcolor{lightgray}1	&	\cellcolor{lightgray}1	\\	
$(a_7,b_5)$	&	0	&	\cellcolor{lightgray}1	&	\cellcolor{lightgray}1	&	\cellcolor{lightgray}1	&	\cellcolor{lightgray}1	&	\cellcolor{lightgray}1	&	\cellcolor{lightgray}1	&	0	&	\cellcolor{lightgray}1	&	\cellcolor{lightgray}1	&	0	&	0	&	0	&	0	&	0	&	0	&	0	&	0	&	0	&	0	&	0	&	0	&	0	&	0	&	0	&	0	&	0	&	0	\\	
$(a_1,b_5)$	&	0	&	0	&	0	&	0	&	\cellcolor{lightgray}1	&	\cellcolor{lightgray}1	&	\cellcolor{lightgray}1	&	0	&	\cellcolor{lightgray}1	&	\cellcolor{lightgray}1	&	0	&	0	&	\cellcolor{lightgray}1	&	0	&	0	&	0	&	0	&	0	&	0	&	0	&	0	&	0	&	0	&	0	&	0	&	0	&	0	&	0	\\	\hline	
$(a_6,b_6)$	&	0	&	0	&	\cellcolor{lightgray}1	&	\cellcolor{lightgray}1	&	0	&	\cellcolor{lightgray}1	&	\cellcolor{lightgray}1	&	\cellcolor{lightgray}1	&	0	&	0	&	0	&	0	&	0	&	0	&	0	&	0	&	0	&	0	&	0	&	0	&	0	&	0	&	0	&	0	&	0	&	0	&	0	&	\cellcolor{lightgray}1	\\	
$(a_7,b_6)$	&	0	&	0	&	\cellcolor{lightgray}1	&	\cellcolor{lightgray}1	&	0	&	\cellcolor{lightgray}1	&	\cellcolor{lightgray}1	&	\cellcolor{lightgray}1	&	\cellcolor{lightgray}1	&	\cellcolor{lightgray}1	&	\cellcolor{lightgray}1	&	0	&	0	&	0	&	0	&	0	&	0	&	0	&	0	&	0	&	0	&	0	&	0	&	0	&	0	&	0	&	0	&	0	\\	
$(a_1,b_6)$	&	0	&	0	&	0	&	0	&	0	&	\cellcolor{lightgray}1	&	\cellcolor{lightgray}1	&	\cellcolor{lightgray}1	&	\cellcolor{lightgray}1	&	\cellcolor{lightgray}1	&	\cellcolor{lightgray}1	&	0	&	\cellcolor{lightgray}1	&	\cellcolor{lightgray}1	&	0	&	0	&	0	&	0	&	0	&	0	&	0	&	0	&	0	&	0	&	0	&	0	&	0	&	0	\\	
$(a_2,b_6)$	&	0	&	0	&	0	&	0	&	0	&	0	&	0	&	0	&	\cellcolor{lightgray}1	&	\cellcolor{lightgray}1	&	\cellcolor{lightgray}1	&	0	&	\cellcolor{lightgray}1	&	\cellcolor{lightgray}1	&	0	&	0	&	\cellcolor{lightgray}1	&	0	&	0	&	0	&	0	&	0	&	0	&	0	&	0	&	0	&	0	&	0	\\	\hline	
$(a_7,b_7)$	&	0	&	0	&	0	&	\cellcolor{lightgray}1	&	0	&	0	&	\cellcolor{lightgray}1	&	\cellcolor{lightgray}1	&	0	&	\cellcolor{lightgray}1	&	\cellcolor{lightgray}1	&	\cellcolor{lightgray}1	&	0	&	0	&	0	&	0	&	0	&	0	&	0	&	0	&	0	&	0	&	0	&	0	&	0	&	0	&	0	&	0	\\	
$(a_1,b_7)$	&	0	&	0	&	0	&	0	&	0	&	0	&	\cellcolor{lightgray}1	&	\cellcolor{lightgray}1	&	0	&	\cellcolor{lightgray}1	&	\cellcolor{lightgray}1	&	\cellcolor{lightgray}1	&	\cellcolor{lightgray}1	&	\cellcolor{lightgray}1	&	\cellcolor{lightgray}1	&	0	&	0	&	0	&	0	&	0	&	0	&	0	&	0	&	0	&	0	&	0	&	0	&	0	\\	
$(a_2,b_7)$	&	0	&	0	&	0	&	0	&	0	&	0	&	0	&	0	&	0	&	\cellcolor{lightgray}1	&	\cellcolor{lightgray}1	&	\cellcolor{lightgray}1	&	\cellcolor{lightgray}1	&	\cellcolor{lightgray}1	&	\cellcolor{lightgray}1	&	0	&	\cellcolor{lightgray}1	&	\cellcolor{lightgray}1	&	0	&	0	&	0	&	0	&	0	&	0	&	0	&	0	&	0	&	0	\\	$(a_3,b_7)$	&	0	&	0	&	0	&	0	&	0	&	0	&	0	&	0	&	0	&	0	&	0	&	0	&	\cellcolor{lightgray}1	&	\cellcolor{lightgray}1	&	\cellcolor{lightgray}1	&	0	&	\cellcolor{lightgray}1	&	\cellcolor{lightgray}1	&	0	&	0	&	\cellcolor{lightgray}1	&	0	&	0	&	0	&	0	&	0	&	0	&	0	\\	\hline	
	\end{tabular}}
	\caption{Adjacency matrix $\mathcal{A}_{4}^{3}$}
	\label{prettytable}
	\end{table}
\end{minipage}
	\begin{minipage}[b]{0.43\textwidth}
	\begin{figure}[H]
\centering
\includegraphics[width=\textwidth]{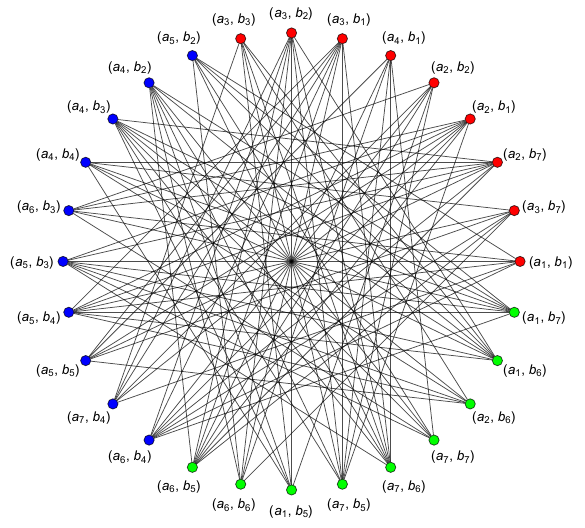}
	 \caption{Graph $\mathbf{G}_4^3\cong\mathcal{BC}_7^4$}
\label{prettygraph}
\end{figure}
\end{minipage} 
\end{minipage}
\end{example}

\begin{theorem}\label{thm1}
Let $\mathcal{A}^k_n=[A_{i,j}]_{1\leq i,j\leq n+k}$, where $n\geq 3$ and  $k\geq 0$.  Then the $(k+1)\times(k+1)$ submatrices $A_{i,j}$ are as follows:
\begin{enumerate}
\item\label{thm1part1} If $i=j$, then $A_{i,j}= \mathbbold{0}_{(k+1)\times (k+1)}$.
\item\label{thm1part2} If $i\neq j$, then $A_{i,j}=[m_{u,v}]$, 
where
\[m_{u,v}=\begin{cases}1&  \mbox{for }u\in\{i,i+1,\ldots,i+k\}\cap\{j+k+1,j+k+2,\ldots,j-1\}\\&
\mbox{ and }v\in\{j, j+1,\ldots,j+k \}\cap\{i+k+1,i+k+2,\ldots,i-1\}\\0&otherwise.\end{cases}\]
\end{enumerate}
  \end{theorem}
  
\begin{proof}
Recall that  $m_{u,v}=1$ if $(b_i,a_u)$ and $(b_j,a_v)$ form a strict alternating cycle, otherwise  $m_{u,v}=0$. By definition, the critical pairs $(a,b)$ and $(a',b')$
  form a strict alternating cycle if and only if 
 (1) $b||a$;
(2)~$a<b'$;
(3) $b'||a'$; and
(4) $a'<b$.
By the definition of $\Skn$,
condition~(1) implies $u\in\{i,i+1,\ldots,i+k\}$;
condition~(2) implies $u\in\{j+k+1,j+k+2,\ldots,j-1\}$;
  condition~(3) implies $v\in\{j,j+1,\ldots,j+k\}$; and
  condition~(4) implies $v\in\{i+k+1,i+k+2,\ldots,i-1\}$.
Thus $m_{u,v}=1$ whenever the preceding 
statements hold simultaneously and otherwise $m_{u,v}=0$ .\\

\noindent{\textbf{Case~1:}}
Assume that $i=j$. 
We claim that $A_{i,i}=\mathbbold{0}_{(k+1)\times(k+1)}$. 
Suppose to the contrary that there exists $u\in[i, i+1,\ldots, k+i]_{\!\!\mod{(n+k)}}$ and $v\in[j, j+1,\ldots, k+j]_{\!\!\mod{(n+k)}}$ such that $m_{u,v}=1$.  
Then $\{ (b_i,a_u),(b_i,a_v)\}$ forms a strict alternating cycle.  Condition~(2) implies that $a_u<b_i $.  This  contradicts  condition~(1). Therefore $A_{i,i}= \mathbbold{0}_{(k+1)\times(k+1)}$. \\

\noindent{\textbf{Case 2:}} If $i\neq j$, then the preceding implications yield the desired result.
  \end{proof}
  
We now state our first result in connection with block circulant matrices.

\begin{theorem}\label{blockcirc} If $n\geq 3$ and $k\geq 0$, then  $\mathcal{A}^k_n$ is a $(k+1)$-block circulant matrix.
\end{theorem}

\begin{proof}
We  show that $A_{i,j}=A_{i+1,j+1}$ for all $1\leq i,j< n+k$. Since $\mathcal{A}_n^k$ is symmetric, we restrict our attention to the case where $i \leq j$.
We proceed by showing that the entries of $A_{i,j}$ are pointwise identical to the entries of $A_{i+1,j+1}$.

Theorem \ref{thm1} states that the nonzero entries of $A_{i,j}$ occur when
 \begin{align*}
 u&\in\{i,i+1,i+2,\ldots,i+k\}\cap\{j+k+1,j+k+2,\ldots,j-1\} \text{ and}\\
 v&\in\{j,j+1,j+2,\ldots,j+k\}\cap\{i+k+1,i+k+2,\ldots,i-1\}, 
 \end{align*}
while the nonzero entries of $A_{i+1,j+1}$ occur when
\begin{align*}
 u'&\in\{i+1,i+2,i+3,\ldots,i+1+k\}\cap\{j+k+2,j+k+3,\ldots,j\} \text{ and}\\
 v'&\in\{j+1,j+2,j+3,\ldots,j+1+k\}\cap\{i+k+2,i+k+3,\ldots,i\}.
 \end{align*}
 
 Fix $i\leq u\leq i+k$ and $j\leq v\leq j+k$. Then the $m_{u,v}$ entry in $A_{i,j}$ corresponds to the $m_{u+1,v+1}$ entry in $A_{i+1,j+1}$.  
Suppose that $m_{u,v}=1$ in $A_{i,j}$. Then  $u\in\{i,i+1,i+2,\ldots,i+k\}\cap\{j+k+1,j+k+2,\ldots,j-1\}$ and $v\in\{j,j+1,j+2,\ldots,j+k\}\cap\{i+k+1,i+k+2,\ldots,i-1\}$.
Hence
\begin{align*}
 u+1&\in\{i+1,i+2,i+3,\ldots,i+1+k\}\cap\{j+k+2,j+k+3,\ldots,j\} \text{ and}\\
 v+1&\in\{j+1,j+2,j+3,\ldots,j+1+k\}\cap\{i+k+2,i+k+3,\ldots,i\}.
 \end{align*}
This implies that $m_{u+1,v+1}=1$ in $A_{i+1,j+1}$.

To complete the proof, note that if $m_{u,v}=0$ in $A_{i,j}$ then
\begin{align*}
    u&\notin\{i,i+1,i+2,\ldots,i+k\}\cap\{j+k+1,j+k+2,\ldots,j-1\} \text{ and}\\
    v&\notin\{j,j+1,j+2,\ldots,j+k\}\cap\{i+k+1,i+k+2,\ldots,i-1\}.
\end{align*}
\noindent
Hence
\begin{align*}
 u+1&\notin\{i+1,i+2,i+3,\ldots,i+1+k\}\cap\{j+k+2,j+k+3,\ldots,j\} \text{ and}\\
 v+1&\notin\{j+1,j+2,j+3,\ldots,j+1+k\}\cap\{i+k+2,i+k+3,\ldots,i\}.
 \end{align*}
This implies that $m_{u+1,v+1}=0$ in  $A_{i+1,j+1}$. Therefore $A_{i,j}=A_{i+1,j+1}$ for any $1\leq i,j\leq n+k$.
\end{proof}
Having shown that the matrices $\mathcal{A}_n^k$ are $(k+1)$-block circulant, when describing $\mathcal{A}_n^k$ we  determine only the generating block row consisting of the $(k+1)\times(k+1)$ matrices $A_{1,1}, A_{1,2}, \ldots, A_{1,n+k}$. 
The following fact will be used throughout the proofs of Propositions \ref{blockcirculant1} and \ref{blockcirculant2}. 

\begin{fact}\label{fact}
Let $X$ and $Y$ be ordered sets with consecutive elements written in increasing order (cyclically). If $Z \subseteq X$ and $Z \subseteq Y$ where no $\alpha \geq \max(Z)+1$ is in both $X$ and $Y$ and no $\beta \leq \min(Z)-1$ is in both $X$ and $Y$, then $Z=X \cap Y$.
\end{fact}

\begin{proposition}\label{blockcirculant1}
Assume $n-1\geq k+1$ and let $1\leq j\leq n+k$.  Then
\begin{align*}
A_{1,j}=&\begin{cases}\mathbbold{0}_{(k+1)\times(k+1)}&\mbox{if $j=1$}\\
\mathbbold{1}^{j-1}&\mbox{if $2\leq j\leq k+1$}\\
\mathbbold{1}_{(k+1)\times(k+1)}&\mbox{if $k+2\leq j\leq n$}\\
_{n+k+1-j}\mathbbold{1}&\mbox{if $n+1\leq j\leq n+k$}.
\end{cases}
\end{align*}
\end{proposition}

\begin{proof}
Theorem~\ref{thm1} states that $A_{1,1}=\mathbbold{0}_{(k+1)\times(k+1)}$.
Next assume that $2\leq j\leq k+1$.  We show that $A_{1,j}=\mathbbold{1}^{j-1}$.  Note that $m_{u,v}=1$ in $A_{1,j}$ if and only if $u\in\{1,2,\ldots,k+1\}\cap\{j+k+1,j+k+2,\ldots,j-1\}$ and $v\in\{j, j+1,\ldots,j+k \}\cap\{k+2,k+3,\ldots,n+k\}$, where each set contains $k+1$ elements listed  (cyclically) in increasing order.

Assume next that $2\leq j\leq k+1$. Then the containments $\{1,2,\ldots, j-1\}\subseteq\{1,2,\ldots,k+1\}$ and $\{1,2,\ldots, j-1\}\subseteq \{j+k+1,j+k+2,\ldots,j-1\}$ hold.
 The latter containment is due to the fact that $\{1, 2, \ldots, j-1\}$ is a set with $j-1 \leq k+1$ elements.
By applying Fact~\ref{fact}, this yields 
$\{1,2,\ldots,k+1\}\cap\{j+k+1,j+k+2,\ldots,j-1\}=\{1,2,\ldots,j-1\}$.
Similarly,
$\{k+2, \ldots, j+k\}\subseteq\{j,\ldots,j+k\}$and $\{k+2,\ldots,j+k\}\subseteq\{k+2,\ldots,n+k\}$ since $j \leq k+1 \leq n-1$, and so $j+k < n+k$. Hence by  Fact~\ref{fact} we have $\{j, j+1,\ldots,j+k \}\cap\{k+2,k+3,\ldots,n+k\}=\{k+2,\ldots,j+k\}$.
Thus, $m_{u,v}=1$ only when $u\in\{1,2,\ldots,j-1\}$ and $v\in\{k+2,\ldots,j+k\}$.  Therefore ${A}_{1,j}=\mathbbold{1}^{j-1}$.

If $k+2\leq j\leq n$, then we show $A_{1,j}=\mathbbold{1}_{(k+1)\times(k+1)}$ by proving that $m_{u,v}=1$ for all possible values of $u$ and $v$.
This follows from  Fact~\ref{fact} and the containments
\begin{align}
\{j,j+1,j+2,\ldots,j+k\}&\subseteq\{k+2,k+3,\ldots,n+k\}\label{subset1} \text{ and}\\
\{1,2,3,\ldots,k+1\}&\subseteq\{j+k+1,j+k+2,\ldots,j+n+k-1\}.\label{subset2}
\end{align}
Containment \eqref{subset1} follows from $k+2\leq j$ and $j+k\leq n+k$. Containment \eqref{subset2} follows from $j \leq n$, which gives $j+k+1 \leq n+k+1 \equiv 1 \!\! \mod{(n+k)}$, and $k+2 \leq j$, which gives $k+1 \leq j-1 \equiv j+n+k-1 \!\! \mod{(n+k)}$. 

If $n+1 \leq j \leq n+k$, we claim $A_{1,j}=\,  _{n+k+1-j}\mathbbold{1}$. It suffices to show that $m_{u,v}=1$ whenever $u\in\{j-n+1,j-n,\ldots, k+1\}$ and $v\in\{j,j+1,\ldots,n+k\}$, and $m_{u,v}=0$ otherwise. Theorem \ref{thm1} states
\begin{align*}
u&\in\{1,2,\ldots,k+1\}\cap\{j+k+1,j+k+2,\ldots,j-1\}\\
v&\in\{j,j+1,\ldots,j+k\}\cap\{k+2,k+3,\ldots,n+k\}.
\end{align*}
Hence it suffices to show that \begin{align}
\{j-n+1,j-n,\ldots, k+1\}&=\{1,2,\ldots,k+1\}\cap\{j+k+1,j+k+2,\ldots,j-1\}\label{a3}\\
\{j,j+1,\ldots,n+k\}&=\{j,j+1,\ldots,j+k\}\cap\{k+2,k+3,\ldots,n+k\}.\label{a4}
\end{align}
Equation \eqref{a3} follows from Fact~\ref{fact} and the inequalities: $n+1\leq j$, which gives $1<j-n+1$, and both $k+1\leq n-1$ and $n+1\leq j$, which gives $k+1<n \leq j-1$. 
Equation \eqref{a4} follows from Fact~\ref{fact} and the inequalities: $n+1\leq j$, which gives $n+k<j+k$, and from both $k+1\leq n-1$ and $n+1\leq j$, which gives $k+2\leq n<j$.
\end{proof}

\begin{proposition}\label{blockcirculant2}
Assume $n-1< k+1$ and let $1\leq j\leq n+k$.  Then

\begin{align*}
A_{1,j}=&
\begin{cases}
\mathbbold{0}_{(k+1)\times(k+1)} &\mbox{if $j=1$}\\
\mathbbold{1}^{j-1} &\mbox{if $2\leq j\leq n$}\\
\begin{bmatrix}0&0&\mathbbold{0}_{(j-n)\times(j-n)}\\0&\mathbbold{1}_{({n-1})\times({n-1})}&0\\
\mathbbold{0}_{({k+2-j})\times({k+2-j})}&0&0
\end{bmatrix}&\mbox{if $n+1\leq j\leq k+1$}\\
_{n+k+1-j}\mathbbold{1} &\mbox{if $k+2\leq j\leq n+k$}.
\end{cases}
\end{align*}
\end{proposition}

\begin{proof}
The equality $A_{1,1}=\mathbbold{0}_{(k+1)\times(k+1)}$ follows from Theorem~\ref{thm1}.   If 
$2\leq j\leq n$, then
 $m_{u,v}=1$ in $A_{1,j}$  when
$u\in\{1,2,\ldots,k+1\}\cap\{j+k+1,j+k+2,\ldots,j-1\}$ and $v\in\{j, j+1,\ldots,j+k \}\cap\{k+2,k+3,\ldots,n+k\}$. Using the inequalities $n-1<k+1$ and $2\leq j\leq n$, we show that 
\begin{align}
\{1,2,\ldots,j-1\}&=\{1,2,\ldots,k+1\}\cap\{j+k+1,j+k+2,\ldots,j-1\}\label{b1}\\
\{k+2,k+3,\dots,j+k\}&=\{j, j+1,\ldots,j+k \}\cap\{k+2,k+3,\ldots,n+k\}.\label{b2}
\end{align}

Equation~\eqref{b1} follows from Fact~\ref{fact} and the observation that 
$2\leq j\leq n$ and $n-1<k+1$ imply
$j-1<k+1$. Equation \eqref{b2} follows from Fact~\ref{fact} and the observation that $2\leq j\leq n$ and $n-1<k+1$ imply $j\leq k+2$.

Next assume  $n+1\leq j\leq k+1$.  Using the inequalities $n-1<k+1$ and $n+1\leq j\leq k+1$, we show 
\begin{align}
\{j-n+1,j-n+2\ldots, j-1\}&=\{1,2,\ldots,k+1\}\cap\{j+k+1,j+k+2,\ldots,j-1\}\label{d1}\\
\{k+2,k+3,\ldots, n+k\}&=\{j, j+1,\ldots,j+k \}\cap\{k+2,k+3,\ldots,n+k\}.\label{d2}
\end{align}
\noindent
Note that $j-n+1\equiv j+k+1 \!\!\mod(n+k)$.  Since $j\leq k+1$ it follows that $j-1< k+1$.  It is clear that $1<j+k+1$.  Therefore 
$\{j+k+1,j+k+2,\ldots, j-1\}\subseteq\{1,2,\ldots,k+1\}$.  By containment and the fact that these sets have the same cardinality, they are equal.
Since $n+1\leq j$, then $n+k<j+k$.  By assumption, $j<k+2$, therefore $\{k+2,\ldots,n+k\}\subseteq\{j,\ldots,j+k\}$.
By containment and the fact that these sets have the same cardinality, they are equal. Thus $m_{u,v}=1$ if and only if $u\in\{j-n+1,\ldots, j-1\}$ and $v\in\{k+2,\ldots, n+k\}$.  

Lastly assume $k+2\leq j\leq n+k$.  Using the inequalities $n-1<k+1$ and $k+2\leq j\leq n+k$, we show 
\begin{align}
\{j-n+1,j-n,\ldots, k+1\}&=\{1,2,\ldots,k+1\}\cap\{j+k+1,j+k+2,\ldots,j-1\}\label{c3}\\
\{j,j+1,\ldots,n+k\}&=\{j,j+1,\ldots,j+k\}\cap\{k+2,k+3,\ldots,n+k\}.\label{c4}
\end{align}
Equation \eqref{c3} follows from Fact~\ref{fact} and the inequalities: $j-n+1\equiv j+k+1\!\! \mod (n+k)>1$, and $k+2\leq j$, the latter implying $k+1\leq j-1$. Equation \eqref{c4} follows from Fact~\ref{fact} and the inequalities $k+2\leq j$ and $n-1<k+1$, which imply $n+k\leq j+k$.
\end{proof}

\begin{theorem}\label{done!}
Let $t\geq 1$ and  $s\geq t+2$. Then any $t$-block circulant matrix $\mathcal{BC}_s^t$ with generating block row as described in Section \ref{sec1} is $\mathcal{A}_{s-t+1}^{t-1}$.
\end{theorem}

The proof of Theorem \ref{done!} follows from the definition of the generating block row of $\mathcal{BC}_s^t$ and from Theorems~\ref{blockcirc}, Proposition \ref{blockcirculant1} and Proposition \ref{blockcirculant2}.

\section{Closing remarks}\label{sec3}

In this paper, we demonstrated a canonical association between a family 
of $m$-block circulant graphs and the classical family of posets known as generalized crowns. It is well-known that the chromatic number of the incomparability graphs is bounded above by the dimension of their associated posets, and we conjecture that this bound is tight. Yet, little has been done in studying the properties of graphs that arise in this way. Naturally, it is of interest to find more families of posets with this property, and if they exist, determine what other attributes their graphs possess. Conversely, the richness in theory and applications surrounding circulant graphs beg their generalization to $m$-block circulant graphs and, specifically, to further investigate the graphs $\mathcal{BC}_s^t$, which are the objects of this study.
 
\section*{Acknowledgments}
The authors acknowledge Charity Bankhead for her initial computational contributions to this project, and thank William T. Trotter for many helpful and enlightening conversations during the completion of this manuscript. As this project began during the 2013 SACNAS National Conference, the authors thank the society for providing an enriching environment that encouraged this scientific collaboration. The authors extend their gratitude to the Center for Leadership and Diversity in STEM at the United States Military Academy, the National Science Foundation Division of Mathematical Sciences (DMS-1045082), and the AMS-Simons Travel Grant for travel support.

\bibliographystyle{amsplain}

\end{document}